\documentclass[11pt]{article}
\usepackage{amsfonts,amssymb,amsmath,pdfsync,color,bbm}
\hsize \textwidth
\advance \hsize by -\marginparwidth
\evensidemargin \oddsidemargin
\usepackage{amssymb}
\usepackage{amsthm}
\usepackage[shortlabels]{enumitem}
\usepackage[normalem]{ulem}
\usepackage{cancel}
\usepackage{scalerel}

\usepackage{mathrsfs}

\advance\hoffset by 5mm
\newcommand{\pdr}[2]{\frac{\partial{#1}}{\partial{#2}}}

\newcommand{\Rm}{{\mathbb R}}

\newcommand{\commentout}[1]{}
\newcommand{\E}{{\mathbb E}}

\newcommand{\Pm}{{\mathbb P}}

\usepackage[colorlinks=true,pdfpagemode=UseNone,urlcolor=blue,linkcolor=blue,citecolor=blue,breaklinks=true]{hyperref}
\usepackage{cleveref}

\newtheorem{thm}{Theorem}[section]

\newtheorem{prop}[thm]{Proposition}

\newcommand{\farc}{\frac}
\newcommand{\be}{\begin{equation}}
\newcommand{\ee}{\end{equation}}

\newcommand{\bal}{\begin{aligned}}
\newcommand{\enbal}{\end{aligned}}

\newcommand{\one}{{\mathbbm{1}}}

\newcommand{\cnk}[2]{\begin{pmatrix}{#1}\cr {#2}\cr\end{pmatrix}}

\newcommand{\vot}{\hbox{Vote}}
\newcommand{\cT}{{\mathcal T}}

\newcommand{\bhat}[1]{\expandafter\hat#1}

\numberwithin{equation}{section}

\begin{document}

\author{Jing An\footnote{Department of Mathematics, Duke University, Durham, NC 27708, USA;
jing.an@duke.edu}  
\and Christopher Henderson\footnote{Department of Mathematics, University of Arizona, 
Tucson, AZ 85721, USA;
ckhenderson@math.arizona.edu} \and 
Lenya Ryzhik\footnote{Department of Mathematics, Stanford University, Stanford, CA 94305, USA;
ryzhik@stanford.edu}}

\title{Voting models and semilinear parabolic equations}

\maketitle

\begin{abstract} We present probabilistic interpretations of solutions to semi-linear parabolic equations
with polynomial nonlinearities in terms of the voting models on the genealogical trees of branching
Brownian motion (BBM).  These extend the connection between the Fisher-KPP equation and BBM discovered
by McKean in~\cite{McK}.  In particular, we present ``random outcome'' and ``random threshold'' voting models that yield any polynomial nonlinearity $f$ satisfying $f(0)=f(1)=0$ 
and a ``recursive up the tree'' model that allows to go beyond this restriction on $f$.  
We compute a few examples of particular interest; for example, we obtain 
a curious interpretation of the heat equation in terms of a nontrivial voting model and  
a ``group-based'' voting rule that leads to a probabilistic view of the pushmi-pullyu transition for a class of nonlinearities introduced by Ebert and van Saarloos.
 \end{abstract}

\section{Introduction}

The Fisher-KPP parabolic equation 
\be\label{aug2616}
u_t=\Delta u+u-u^2,
\ee
was introduced as a minimal model for gene mutations in the pioneering works by Fisher~\cite{Fisher} and Kolmogorov, Petrovskii and Piskunov~\cite{kpp} in their seminal
papers that both appeared in 1937. More general semi-linear parabolic equations of the form 
\be\label{aug2618}
u_t=\Delta u+f(u),
\ee
have been used in areas as diverse as biology, combustion, and differential geometry, as well as many others.  Quite often they are considered
in a context of a diffusive connection between two equilibria states that may be assumed, without loss of generality, to be $0$ and $1$, so that
\be\label{aug2620}
f(0)=f(1)=0.
\ee

In 1975, McKean~\cite{McK} discovered an elegant interpretation of the solutions to the original Fisher-KPP equation 
(\ref{aug2616}) in terms of branching Brownian motion (BBM).  While we recall this in greater detail in \Cref{sec:mckean}, in order to illustrate this connection,
let us just mention that  if $M(t)$ is the running maximum of the binary BBM that starts at a location $x=0$ at the time
$t=0$, then the probability distribution
function 
\be\label{aug2621}
u(t,x)=\Pm(M(t)>x),
\ee
is the solution to (\ref{aug2616}) with the initial condition $u(0,x)=\one(x<0)$. McKean's interpretation was used both to understand the
long time behavior of the solutions to (\ref{aug2616}) via a study of the fine properties of BBM, as in~\cite{Bovier,Bramson1,Bramson2}, and to
understand the properties of BBM in terms of the solutions to the Fisher-KPP equation~\cite{BD1,BD2,MRR}.

McKean's interpretation directly extends to solutions to equations (\ref{aug2618}) for a larger set of nonlinearities than $f(u)=u-u^2$, known as the
McKean class. 
However, while McKean's interpretation is extremely beautiful,  
the McKean class of nonlinearities is still rather limited, and equations (\ref{aug2618})  with a McKean type nonlinearity $f(u)$
have some very special properties. For example, their solutions are all pulled, in the sense that their spreading is governed by the behavior in the regions
where $u(t,x)$ is very small.
This is also discussed further in Section~\ref{sec:mckean} below.  

The standard example of a semi-linear parabolic equation that does not admit a McKean interpretation
is the Allen-Cahn equation
\be\label{aug2622}
u_t=\Delta u+u(1-u)\Big(u-\farc{1}{2}\Big).
\ee
It appears both as a basic model of a diffusive connection between two stable states $u=0$ and~$u=1$ in physics and biology, as well
as in the analysis of mean curvature flow, and it is closely connected to the De Giorgi conjecture. 
It was believed for a long time that solutions to the Allen-Cahn equation do not have any probabilistic interpretation; however,
recently, Etheridge, Freeman, and Penington discovered in~\cite{EFP} that they can be represented using 
ternary branching Brownian motion. Their key idea, recalled
in Section~\ref{sec:efp},  is to describe the solution not in terms of just the locations of the BBM particles, as in the McKean argument, 
but to use a voting procedure on the genealogical
tree of the~BBM.

In this paper, we explore their ideas further and show that solutions to any semi-linear parabolic equation (\ref{aug2618}) with a polynomial nonlinearity
can be described
in terms of recursive procedures on the genealogical trees of the BBM. For nonlinearities that satisfy (\ref{aug2620}), we describe two randomized voting
models: a random outcome model and a random threshold model proposed in~\cite{KRZ}.  In the former, the vote of the parent particle is influenced but 
not deterministically determined by the votes of the children. In the latter, the vote of the parent is determined by the votes of its children. However, 
the minimal number of children who vote $1$ that is required for the parent to vote $1$, is random at each vertex of the genealogical tree.
Somewhat amusingly, these models also give what seems to be a new probabilistic interpretation for the solutions to the standard heat equation,
in terms for a BBM rather than as its classical description as the Kolmogorov equation for the standard Brownian motion.  
Finally,  we introduce a recursive procedure up the genealogical tree of the BBM that gives a
probabilistic interpretation for solutions to equations (\ref{aug2618}) with polynomial nonlinearities that need not satisfy assumption~(\ref{aug2620}). 
We should mention that the random outcome voting model turns out to be a special case 
of the collection of models considered in the impressive thesis of O'Dowd~\cite{ODowd}
who obtained results similar to ours
in that case. His thesis, unfortunately, is not widely available. 

Let us also recall that BBM is but one example of a log-correlated random process that includes objects as diverse as eigenvalues 
of random matrices~\cite{ABB,FHK1,FHK2,CMN,Lam-Paq,PZ},
extremal values of the Gaussian free field~\cite{DZ,Zeit-Notes}, and local maxima of the zeta function on large intervals on the critical
line~\cite{AOR,ABBRS,ABR,Naj}, 
that all have a natural tree structure. It is known that
their extremal values all belong to the same universality class. 
The results of this paper give a tool to study various other, non-extremal, observables of the BBM trees. One may hope that this allows one to understand the 
extent of the aforementioned universality of extremal values of the log-correlated processes, and how far it can be extended away from the purely extremal values. 
In the last section of the paper we describe a voting model for BBM with an odd number of offspring
that includes such an extremal-to-bulk transition. It is tempting to conjecture that the nature
of this transition
is generic for the observables of the log-correlated processes.

The paper is organized as follows. In Section~\ref{sec:mckean} we recall the McKean interpretation of the solutions to the Fisher-KPP
equation (\ref{aug2616}) in terms of the BBM. The main results of this paper and the description of the voting schemes are presented in
Section~\ref{s.voting}. Some concrete examples are discussed in Section~\ref{sec:examples}. 

Let us mention that all Brownian motions in this paper run with the 
diffusivity $\sigma=\sqrt{2}$, so that the pre-factor in front of the Laplacian in (\ref{aug2616}) is $1$ and not $1/2$.  

\noindent{\bf Acknowledgment.} We are greatly indebted to Sarah Penington for 
numerous illuminating discussions, and to Alison Etheridge for pointing out O'Dowd's thesis
to us.
CH was supported by NSF grants DMS-2003110 and DMS-2204615. 
LR was supported by NSF grants DMS-1910023 and DMS-2205497, and by ONR grant N00014-22-1-2174.

\section{McKean nonlinearities and branching Brownian motion}\label{sec:mckean}

\subsection{The McKean interpretation}

In this section, we recall McKean's interpretation~\cite{McK} for solutions to a class of semilinear parabolic equations of the form
\be\label{aug402}
u_t=\Delta u+f(u),
\ee
in terms of BBM. 
Let us briefly recall that BBM is defined as follows. A single particle starts  at a position~$x\in\Rm^d$ at $t=0$ and performs a standard Brownian motion
$B_0(t)$. The particle
carries an exponential clock that rings at a random time $\tau$, with 
\be\label{aug404}
	\qquad
	\Pm(\tau>t)=e^{-\beta t}
	\qquad\text{ for some } \beta>0.
\ee
At the time $\tau$ the particle splits into a random number $N\ge 2$ of ``offspring'' particles that we refer to as the ``children''.
The original
particle is sometimes called the ``parent''.
The probabilities 
\be\label{aug406}
p_k=\Pm(N=k),
	\quad k\ge 2,
\ee 
are fixed, and
\be\label{aug2302}
\sum_{k=2}^\infty p_k=1.
\ee
The original particle is removed at the branching event. The $k$ children 
perform independent standard Brownian motions for~$t>\tau$, and all of them starting at the position $B_0(\tau)$ of the branching event. 
Each of the children carries its own exponential clock, and when
the corresponding clock rings, that particle splits into a random number of particles, with the same probabilities $p_k$, and the process continues. 
Thus, at each time $t>0$ we have a collection of particles~$x_1(t),\dots,x_{N_t}(t)$.  Detailed properties of BBM are discussed in~\cite{JB,Bovier}. 

An insightful observation of McKean 
is that if we fix a function~$g\in C_b(\Rm^d)$, then the function
\be\label{aug408}
v(t,x)=\E\Big(\prod_{m=1}^{N_t}g(X_m(t))\Big)
\ee
is the solution to the initial value problem
\be\label{aug410}
\bal
&v_t=\Delta v+\beta F(v),\\
&v(0,x)=g(x).
\enbal
\ee
The nonlinearity $F(v)$ in (\ref{aug410}) is explicit:
\be\label{aug412}
F(v)=\sum_{k=2}^\infty p_k(v^k-v).
\ee
In practice, we restrict to the case where $p_k = 0$ for all but finitely many $k$.  We denote by $N$ the largest possible number of children; that is, the largest $k$ such that $p_k > 0$.   
A simple way to see that~(\ref{aug410}) holds is by writing the renewal equation for $v(t,x)$ using its definition~(\ref{aug408})
and the fact that the trajectories are Markovian, as well as the independence of the children from each other:
\be\label{aug414}
\bal
v(t,x)&=\E(g(x+B_t))\Pm(\tau>t)+ \int_0^t \E\Big(\sum_{k=2}^N 
p_kv^k(t-s,x+B_s)\Big)\Pm(\tau\in ds)\\
&=\E(g(x+B_t))e^{-\beta t}+ \int_0^t \E\Big(\sum_{k=2}^N p_kv^k(t-s,x+B_s)\Big)\beta e^{-\beta s}ds\\
&=e^{(\Delta-\beta)t}g(x)+\beta  \int_0^t \Big(\sum_{k=2}^N 
p_ke^{(\Delta-\beta)s} v^k(t-s,\cdot)(x)\Big) ds.
\enbal
\ee
This is the Duhamel formulation of (\ref{aug410}). The first term in the right side of (\ref{aug414}) accounts for the event that there was no branching
until the time $t$, so that there is only one particle present at $t$. The second term accounts for the event that the first branching happened
in a time interval~$[s,s+ds]$ with $0<s\le t$, and uses the product structure of the definition (\ref{aug408}) of the function $u(t,x)$ and
independence of the offspring particles.  

In the PDE literature, it is more common to consider the
function $u=1-v$ rather than $v$ itself. It satisfies an equation of the form
\be\label{mar1630}
\bal
&u_t=\Delta u+f(u),\\
&u(0,x)=1-g(x).
\enbal
\ee
with the nonlinearity
\be\label{mar1631}
f(u)=\beta\Big(1-u -\sum_{k=2}^N p_k(1-u)^k\Big). 
\ee
The nonlinearities of the form (\ref{mar1631}) are often called the McKean nonlinearities.

\subsection{The Fisher-KPP and McKean nonlinearities}\label{sec:fkpp-prop}

Let us discuss some of the  properties of the McKean nonlinearities.
First, because $p_k$ are probabilities, so that (\ref{aug2302}) holds, the definition
(\ref{mar1631}) of~$f(u)$ immediately implies that 
\be\label{oct518}
	f(0)=f(1)=0
	\quad\text{ and }\quad
	f(u)>0
	~~ \hbox{ for $0<u<1$.}
\ee
Second, it follows from (\ref{mar1631}) that
$f(u)$ is concave. In particular, it necessarily 
satisfies the ``Fisher-KPP condition''
\be\label{oct520}
f(u)\le f'(0)u,~~\hbox{ for all $u\in(0,1)$.}
\ee
This means that the steady state $u=1$ is stable and $u=0$ is an unstable steady state for (\ref{aug2616}). Thus,  solutions to (\ref{aug2616}) 
with a non-negative compactly supported initial condition $u(0,x)=g(x)$ are spreading in the sense that the state $u\approx 1$ is invading the
region where $u\approx 0$. The front location for the solutions to (\ref{aug2616}) with a Fisher-KPP type nonlinearity $f(u)$ 
has the asymptotics
\be\label{aug2624}
	X(t)
		=2\sqrt{f'(0)}t-\farc{3}{2\sqrt{f'(0)}}\log t+x_0
	\quad \hbox{ as } t\to +\infty.
\ee
This was first established by Bramson in~\cite{Bramson1,Bramson2} and later studied in~\cite{berestycki2018new,Graham,HNRR,Lau,NRR1,Roberts,Uchiyama},
among others. The coefficient $3/2$ in (\ref{aug2624}) is a universal signature of the extremals of 
log-correlated processes that we have discussed in the introduction.
In PDE terms, it is present in the ``pulled'' regime that was also mentioned 
in the introduction~\cite{AHR-RD,AS2,Ebert-vanSaarlos,garnier2012inside,Zlatos-lin}. 
However, the Bramson asymptotics~(\ref{aug2624}) for the
solutions to~(\ref{aug2616}) extends beyond the McKean or Fisher-KPP nonlinearities~\cite{AHR-RD,AS2}. Thus, it is natural to expect that for other log-correlated fields
this asymptotics  holds within some range that is away from the ``purely extremal'' observables as well. 


\subsection{Not all Fisher-KPP nonlinearities are of the McKean type}\label{sec:mck-notkpp}

As we have noted, the McKean nonlinearities belong to the class of
Fisher-KPP nonlinearities, and solutions to parabolic equations (\ref{aug402})
with a FKPP type $f(u)$ enjoy many special properties. However, while
the McKean nonlinearities lie in the FKPP  class, they form a very special
sub-class of that set which excludes many natural examples. 
  
The McKean nonlinearities in (\ref{mar1631}) can be written in the form 
\be\label{mar708}
f(u)=\lambda (u-A(u)),
\ee
with
\be
\lambda
	=\beta\Big(\sum_{k=2}^N k p_k -1\Big)
	=f'(0)>0,
\ee
and  the function $A(u)$ defined by
\be
\lambda A(u)=\lambda u-f(u)
	=\beta\sum_{k=2}^N p_k \Big((1-u)^k-1+ku\Big).
\ee
One can immediately check that $A(u)$ is non-negative and convex on $[0,1]$, and 
\be\label{mar718}
	A(0)=0,
	\qquad A(1)=1,
	\qquad A'(0)=0.
\ee 
It follows, in particular, that  the function $A(u)$ is increasing on $[0,1]$.  

It is convenient for us to write $f(u)$ in the form
\be\label{mar702}
f(u)=\lambda u(1-\alpha(u)),
\ee
with 
\be\label{mar704}
\alpha(u)=\farc{A(u)}{u}= \farc{\beta}{\lambda}\sum_{k=2}^N p_k\alpha_k(u).
\ee
The coefficients $\alpha_k(u)$ are 
\be\label{mar706}
\alpha_k(u)= \farc{(1-u)^k-1+ku}{u}=k-\farc{1-(1-u)^k}{u}=k-(1+(1-u)+\dots+(1-u)^{k-1}).
\ee
Note that each $\alpha_k(u)$ is increasing and concave on $[0,1]$ with $\alpha_k(0)=0$ and $\alpha_k(1)=k-1$. 
It follows from (\ref{mar702})-(\ref{mar706}) that if $f(u)$ is a McKean nonlinearity then 
\be\label{mar716}
\alpha(0)=0,~~ \alpha(1)=1,
\ee
and $\alpha(u)$ is increasing and concave. 
 
The most standard example of the Fisher-KPP nonlinearity is $f(u)=u-u^2$,
which is in the McKean class (\ref{mar1631}), as can be
seen by setting $p_2=1$ and $p_k=0$ for $k>2$. However, even the original example
\be\label{aug2312}
f(u)=u(1-u)^2=u-2u^2+u^3
\ee
in the
KPP paper~\cite{kpp} is not of that form because the function
\be
A(u)=2u^2-u^3
\ee
is not convex for $u\in(2/3,1)$. Moreover, the functions $f_n(u)=u-u^n$ are of the
FKPP type but the corresponding functions $\alpha_n(u)=u^{n-1}$ are convex and not
concave. This means that they are also not of the McKean type. As a consequence, McKean's
connection between semi-linear parabolic equations and branching Brownian motion does
not cover all polynomial Fisher-KPP nonlinearities.

\section{Voting schemes for semi-linear equations}\label{s.voting}

In this section, we describe a probabilistic interpretation for 
semilinear parabolic equations with polynomial nonlinearities, in terms of voting models
for branching Brownian motion. These results are inspired by the arguments in~\cite{EFP,EP} for the
Allen-Cahn nonlinearity.  We first recall that connection in Section~\ref{sec:efp}. Then, we discuss the random outcome and
random threshold voting models that allow us to give a probabilistic interpretation to solutions to (\ref{aug2616})
with any polynomial nonlinearity~$f(u)$ such that 
\be\label{aug2625}
f(0)=f(1).
\ee
This is done in Sections~\ref{sec:alpha} and \ref{sec:thresh},
with the main results given in Theorems~\ref{thm-aug2402} and~\ref{thm-aug2404}. Section~\ref{sec:alpha} also describes
an interpretation   of the standard heat equation, in terms of a BBM and an unbiased voting model, formulated in Proposition~\ref{thm-aug2604}.  
Finally, in Section~\ref{sec:recursive} we drop the assumption (\ref{aug2625})  and describe a recursive procedure on the genealogical tree
that gives a BBM-interpretation  for the solutions to any equation of the form (\ref{aug2616}) with a polynomial
nonlinearity $f(u)$. The result is described in Theorem~\ref{thm-aug2602}.

 \subsection{The Etheridge-Freeman-Penington model for the Allen-Cahn equation}\label{sec:efp}

An alternative connection between semilinear parabolic equations and branching Brownian motion to McKean's was pointed out in a beautiful 
paper by Etheridge, Freeman, and Pennington~\cite{EFP}.  One of the main points of the approach of~\cite{EFP,EP}
is to consider functionals that depend not just on
the locations of the BBM particles, as was done by McKean in~\cite{McK}, but also on the  
structure of the (random) genealogical tree that results from the branching.  
There is a natural way to associate a random genealogical tree $\cT(t)$ to each realization of
the BBM running on a time interval~$0\le s\le t$. Each vertex of the tree corresponds to a branching event, while each of the edges
coming out of a vertex represents an offspring particle born at that branching event. The root of the tree $\cT(t)$ represents the original particle
that started at the time $s=0$ at a position~$x\in\Rm^d$. We refer to~\cite{EFP,EP}
for a formal definition of~$\cT(t)$.

Before introducing a generalization of their ideas, let us recall the example of~\cite{EFP}.
Consider a ternary branching Brownian motion starting at the time~$t=0$ at a point~$x\in\Rm^d$
-- each branching event produces three children. The process is run
until a time $t>0$, with the BBM particles at the time $t>0$ located
at the positions~$X_1(t),\dots,X_{N_t}(t)$.  
Then, each of the youngest generation particles~$X_j(t)$,~$j=1,\dots,N_t$,  
``votes"~$0$ or $1$, with the probabilities
\be\label{oct616}
	\Pm(V_j=1)=g(X_j(t))
	\quad\text{ and }\quad
	\Pm(V_j=0)=1-g(X_j(t)).
\ee
Here, $g(x)$ is a prescribed function such that $0\le g(x)\le 1$ for all $x\in\Rm$, and $V_j$ is the vote of the particle~$X_j(t)$, $j=1,\dots, N_t$.
This produces the votes of the youngest generation of particles.  
Next, we go back up the ternary branching tree $\cT(t)$, with the rule that each parent
accepts the vote of the majority of its three children. In this way, we obtain the votes of all particles
on the genealogical tree. 

Let $V_{\rm orig}$ be the resulting vote of the original ancestral particle
that started at $t=0$ at the position $x$, and consider the function  
\be\label{oct544}
u(t,x)=\Pm_x(V_{\rm orig}=1).
\ee
We now derive an equation for $u(t,x)$ using a similar approach to  (\ref{aug414}). 
There are exactly two possible ways
in which the original ancestor can vote $1$: either all three of its children voted $1$ or two of them voted $1$ and one
voted $0$. In the latter case, there are three choices of the particle that voted~$0$. If there has been no branching 
before the time $t$ then the only particle present is the original ancestor, and it takes the vote by itself. 
This gives the renewal identity
\be\label{oct545}
\bal
u(t,x)&=\E_x\big[g(B_t)\big]\Pm(\tau_1>t)\\
&+\int_0^t \E_x(u^3(t-s,B_s)+3u^2(t-s,B_s)(1-u(t-s,B_s))\Pm(\tau_1\in ds)\\
&=\E_x[g(B_t)]e^{-\beta t}+\beta\int_0^t \E_x(u^3(t-s,B_s)+3u^2(t-s,B_s)(1-u(t-s,B_s))e^{-\beta s}ds\\
&=e^{(\Delta-\beta)t}g(x)+\beta\int_0^t e^{(\Delta-\beta)s}\big[u^3(t-s,\cdot)+3u^2(1-u)(t-s,\cdot)\big](x)ds.
\enbal
\ee
A simple computation shows that
\be
u^3+3u^2(1-u)-u=3u^2-2u^3-u=u(3u-2u^2-1)=u(1-u)(2u-1).
\ee
We deduce that the function $u(t,x)$ defined by (\ref{oct544}) satisfies the Allen-Cahn equation
\be\label{oct547}
\bal
&\pdr{u}{t}= \Delta u +u(1-u)(2u-1),\\
&u(0,x)=g(x).
\enbal
\ee
This equation is probably the most standard example of a semi-linear parabolic equation that does not have
a McKean connection to BBM and, prior to~\cite{EFP}, was believed to have no probabilistic interpretation,
to the best of our knowledge. 
Note that the nonlinearity 
\be\label{oct827}
f(u)=u(1-u)(2u-1)
\ee
does not satisfy the Fisher-KPP properties we have discussed in Section~\ref{sec:fkpp-prop}. Indeed, it is not even non-negative
for $u\in(0,1)$ but rather changes its sign. Thus, the Allen-Cahn equation~(\ref{oct547}) does not have an interpretation in terms of
a McKean functional. The voting scheme idea of~\cite{EFP} adds a genuinely new aspect here and dramatically broadens the class of equations that have
an interpretation in terms of the BBM. 

\subsection{Random outcome probabilistic voting models} \label{sec:alpha}

The voting procedure of~\cite{EFP} that we have described above is deterministic, in the sense that once the genealogical tree $\cT(t)$ and the votes
of the youngest generation particles $X_1(t),\dots,X_{N_t}(t)$ are fixed, the vote $V_{\rm orig}$ of the original particle is completely determined. 
However, one can also randomize the voting process itself, in at least  two clear ways that we now discuss.  We  consider a general BBM, with the
probabilities $p_k$ to produce $k$ offspring particles at each branching event, with~$2\le k\le N$. As before, we  denote by $\cT(t)$ the genealogical tree
produced by branching events on the time interval~$0\le s\le t$, and by $p_k$ the probability that a parent produces exactly $k$ children
at a given branching event. 

Let us fix a continuous function $g(x)$ such that $0\le g(x)\le 1$ for all $x\in\Rm^d$ and run a branching Brownian motion starting at a position~$x\in\Rm^d$ until a time 
$t>0$. At the time $t$, each of the BBM particles~$X_1(t),\dots,X_{N_t}(t)$ 
votes randomly~$0$ or~$1$, with the probability to vote $1$ given by 
\be\label{mar1702}
\Pm(\hbox{Vote}(X_k(t))=1)=g(X_k(t)),~~\hbox{for each $1\le k\le N_t$.}
\ee
In a difference with~\cite{EFP},  we also fix a collection of probabilities $0\le \alpha_{kn}\le 1$, with~$0\le k\le n$,
and~$n\ge 2$, such that
\be\label{aug2314}
\alpha_{0n}=0,~~\alpha_{nn}=1,~~\hbox{ for all $n\ge 2$.}
\ee
Given the votes of the particles that are present at the time $t$, we propagate the vote up the genealogical tree $\cT(t)$ as follows. If a parent  
particle has $n$ children and $k$ out of its $n$ children voted~$1$, then the parent particle votes $1$ with the probability $\alpha_{kn}$.
That is, the vote of the parent is no longer a deterministic function of the votes of its children. 
Using this rule to go up the tree all the way to the root  produces the vote~$\vot_{\rm orig}$ of the original ancestor particle, and we can, as before, define
\be\label{mar1704}
u(t,x)=\Pm_x(\vot_{\rm orig}=1).
\ee
If there was no branching event until the time $t$, so that $N_t=1$, then the vote of the original particle is $1$ with the probability $g(X_1(t))$. 
Note that while we allow $\alpha_{kn}$ to be different from $0$ and $1$, we do impose (\ref{aug2314}) which says that if all children 
voted unanimously, then the parent accepts the vote of the children. We  refer to the above as a random outcome voting model. 

Similarly to (\ref{oct545}), one  can write a renewal equation for the function $u(t,x)$ defined in (\ref{mar1704}):
\be\label{mar1706}
\bal
&u(t,x)=\E(g(x+B_t))\Pm(\tau>t)\\
&\quad + \int_0^t \E\Big(\sum_{n=2}^N p_n\sum_{k=0}^n\cnk{n}{k}\alpha_{kn}u^k(t-s,x+B_s)(1-u(t-s,x+B_s))^{n-k}\Big)\Pm(\tau\in ds)\\
&=\E(g(x+B_t))e^{-\beta t}\\
&\quad+ \int_0^t  \E\Big(\sum_{n=2}^Np_n\sum_{k=0}^n\cnk{n}{k}\alpha_{kn}u^k(t-s,x+B_s)(1-u(t-s,x+B_s))^{n-k}\Big)\beta e^{-\beta s}ds\\
&=e^{(\Delta-\beta)t}g(x)+\beta  \int_0^t \Big( \sum_{n=2}^N p_n
\sum_{k=0}^n\cnk{n}{k}\alpha_{kn} e^{(\Delta-\beta)s} u^k(t-s,\cdot)(1-u(t-s,\cdot))^{n-k}(x)\Big) ds.
\enbal
\ee
The first term on the right in (\ref{mar1706}) comes from the event that there was no branching until the time
$t$. The second accounts for the first branching event happening at a time $t\in[s,s+ds]$,
with~$0<s<t$. The binomial coefficient counts the number of possibilities to choose the
$k$ children who voted $1$ out of the $n$ children. 
Note that (\ref{mar1706}) is the Duhamel formulation of the initial value problem
\be\label{mar1708}
\bal
&u_t=\Delta u+f(u),\\
&u(0,x)=g(x),
\enbal
\ee
with the nonlinearity
\be\label{mar1710}
\bal
f(u)&=\beta\sum_{n=2}^N p_n\sum_{k=0}^n\cnk{n}{k}\alpha_{kn}u^k(1-u)^{n-k} -\beta u=
\beta\sum_{n=2}^N p_n\Big(\sum_{k=0}^n\cnk{n}{k}\alpha_{kn}u^k(1-u)^{n-k} -u\Big).
\enbal
\ee

As we have seen in the Etheridge-Freeman-Penington example, 
unlike in the McKean interpretation, the nonlinearities produced in this way need not  be of the Fisher-KPP 
type. 
The advantage of the voting models is precisely in providing a
probabilistic interpretation for nonlinear parabolic 
equations not accessible by the McKean formula.

\subsubsection*{The standard heat equation and unbiased voting}
 
We now consider some  concrete examples of parabolic equations coming from
probabilistic voting models, starting with the standard heat equation. 
Let us first note an elementary identity:  for any~$n\ge 1$ we have
\be\label{jan2814}
\sum_{k=0}^n \cnk{n}{k}\farc{k}{n}u^k(1-u)^{n-k}=u.
\ee
To see why (\ref{jan2814}) holds, we re-write the left side   as
\be\label{jan2816}
\bal
\sum_{k=0}^n\farc{n!}{k!(n-k)!}&\farc{k}{n}u^k(1-u)^{n-k}=\sum_{k=1}^n\farc{(n-1)!}{(k-1)!(n-k)!} u^k(1-u)^{n-k}\\
&=\sum_{k=0}^{n-1}\farc{(n-1)!}{k!(n-1-k)!} u^{k+1}(1-u)^{n-1-k}=u\sum_{k=0}^{n-1}\cnk{n-1}{k}u^k(1-u)^{n-1-k}\\
&=u(u+1-u)^{n-1}=u.
\enbal
\ee
 
Using (\ref{jan2814}) in the representation (\ref{mar1710}) for $f(u)=0$, we see that taking  the probabilities 
\be\label{mar1716}
\alpha_{kn}=\farc{k}{n}
\ee
in the above voting scheme leads to the standard heat equation: (\ref{mar1708}) becomes
\be\label{mar1718}
\bal
&u_t=\Delta u,\\
&u(0,x)=g(x).
\enbal
\ee
That is, consider any branching Brownian motion, regardless of the branching probabilities~$p_k$, 
and introduce the voting scheme such that a parent with 
$n$ children, out  of which $k$ voted $1$, votes~$1$ with the ``unbiased'' probability $\alpha_{kn}=k/n$. Then, the
function $u(t,x)$, the probability that the original ancestor particle votes $1$, is the solution to the standard heat equation.
To the best of our knowledge, even this very simple and intuitive probabilistic 
interpretation of the heat equation is new.  Let us summarize this result as follows. 
\begin{prop}\label{thm-aug2604}
Let $g(x)$  be a continuous function that satisfies~$0\le g(x)\le 1$
for all $x\in\Rm^d$. Consider the random outcome voting model with the  voting probabilities
$\alpha_{kn}=k/n$, $0\le k\le n$, for any branching Brownian motion. Then, the function $u(t,x)=\Pm_x(V_{\rm orig}=1)$   
is the solution to the initial value problem
\be\label{aug2626}
\bal
&u_t=\Delta u,~~t>0,~x\in\Rm^d,\\
&u(0,x)=g(x),~~x\in\Rm^d.
\enbal
\ee
\end{prop}

\subsubsection*{Representing general nonlinearities}

Let now $f(u)$ be a polynomial of degree $N$:
\be\label{aug2404}
f(u)=\sum_{k=0}^N f_ku^k,
\ee
that vanishes at $u=0$ and $u=1$:
\be\label{aug2406}
f(0)=f(1)=0.
\ee
Our goal is to find $\beta>0$ and $\alpha_{kN}$, $0\le k\le N$ such that 
\be\label{aug2408}
\alpha_{0N}=0,~~\alpha_{NN}=1,~~0\le\alpha_{kN}\le 1,~~\hbox{ for all $1\le k\le N-1$,}
\ee
and so that representation (\ref{mar1710}) holds for $f(u)$. We  set all $p_n=0$ except for $p_N=1$, so that each branching event
produces exactly $N$ offspring. We look for $\alpha_{kN}$ in the form 
\be\label{aug2410}
\alpha_{kN}=\farc{k}{N}+  \mu_k
	\quad\text{ for any } 1\le k\le N-1,
\ee
with $\mu_k$ to be chosen later. Recalling ~\eqref{mar1710} and~(\ref{jan2814}), we see that we need to have
\be\label{aug2412}
f(u)= \beta \sum_{k=1}^{N-1}\cnk{N}{k}\mu_{k}u^k(1-u)^{N-k}.
\ee
Recall that the Bernstein polynomials
\be\label{aug2414}
B_{k,N}(u)=\cnk{N}{k}u^k(1-u)^{N-k},~~k=0,\dots,N,
\ee
form a basis for the vector space of polynomials of degree at most $N$. Therefore, $f(u)$ has a representation
\be\label{aug2416}
f(u)=\sum_{k=0}^Nb_k[f]B_{k,N}(u),
\ee
with some coefficients $b_k[f]$. Note that
\be\label{aug2418}
b_0[f]=f(0),~~b_N[f]=f(1).
\ee
We deduce from (\ref{aug2406}) and (\ref{aug2418}) that 
\be\label{aug2420}
b_0[f]=b_N[f]=0. 
\ee
Next, comparing (\ref{aug2412}) and (\ref{aug2416}), we see that for (\ref{aug2412}) to hold, we need to have 
\be\label{aug2422}
\mu_k=\farc{b_k[f]}{\beta  },~~\hbox{for all $1\le k\le N-1$.}
\ee
It remains to choose  $\beta>0$ so that $\alpha_{kN}$ given by (\ref{aug2410}) satisfy (\ref{aug2408}). Note that, 
since we have set~$\mu_0=\mu_N=0$, we automatically have $\alpha_{0N}=0$ and $\alpha_{NN}=1$. The rest of the conditions in (\ref{aug2408}) 
translates into  
\be\label{aug2423}
0\le \farc{k}{N}+ \mu_k\le 1,~~\hbox{for all $1\le k\le N-1$.} 
\ee
We see from  (\ref{aug2422}) that this is equivalent to 
\be\label{aug2424}
0\le \farc{k}{N}+\farc{b_k[f]}{\beta}\le 1,~~\hbox{for all $1\le k\le N-1$.} 
\ee
This condition holds as long as we choose $\beta>0$ sufficiently large, so that 
\be\label{aug2425}
\beta\ge N\max_{k}|b_k[f]|.
\ee
Therefore, we have proved the following.
\begin{thm}\label{thm-aug2402}
Let $f(u)$ be a polynomial of degree $N$ such that $f(0)=f(1)=0$. Then, there exists a random outcome voting model representation  in terms of a purely
$N$-ary branching Brownian motion for the 
solution to the initial value problem (\ref{mar1708}) with the initial condition $g(x)$ that is continuous and satisfies~$0\le g(x)\le 1$
for all $x\in\Rm^d$. 
\end{thm}

We note that \Cref{thm-aug2402} has also been observed in~\cite{ODowd}, although with different terminology and presentation.

\subsection{A random threshold   voting model} \label{sec:thresh}

An alternative probabilistic voting model has been recently suggested in~\cite{KRZ}. For simplicity of notation, let us fix the number~$N$ of offspring 
produced at each branching event and consider a purely $N$-ary BBM. Then, at each branching event we choose a number $L\in\{1,\dots,N\}$,
with the probability $\zeta_k=\Pm(L=k)$ so that
\be\label{aug2427}
\sum_{k=1}^N\zeta_k=1.
\ee
Thus, an integer $L(\nu)$ is assigned separately to each vertex $\nu$ of the genealogical tree $\cT(t)$. 
The voting is done as follows. As before, at the time $t$ the youngest generation of particles votes according to~(\ref{mar1702}). 
The difference is in the way the votes are propagated up the genealogical tree. A parent at a vertex $\nu$ votes $1$
if and only if at least $L(\nu)$ of its children voted $1$. We refer to this process as a random threshold voting model.

The same argument as in (\ref{mar1706}) shows that the function
\be\label{aug2429bis}
u(t,x)=\Pm_x(\vot_{\rm orig}=1)
\ee
is a solution to the initial value problem
\be\label{aug2430}
\bal
&u_t=\Delta u+G(u),\\
&u(0,x)=g(x),
\enbal
\ee
with 
the nonlinearity
\be\label{jun2114}
G(u)=\beta\sum_{j=0}^N\zeta_{j}\sum_{k=j}^N\cnk{N}{k}u^k(1-u)^{N-k}-\beta u=
\beta\sum_{k=0}^N \cnk{N}{k}u^k(1-u)^{N-k}\sum_{j=0}^k\zeta_{j}-\beta u.
\ee  
A simple observation is that if we start with the random threshold voting model and set 
\be\label{jun2118}
\alpha_{kN}=\sum_{j=0}^k\zeta_{j},
\ee
in (\ref{mar1710}), then the nonlinearities $f(u)$ in (\ref{mar1710}), coming from the random outcome model with the probabilities $\alpha_{kN}$,
and $G(u)$ in (\ref{jun2114}) are the same. 
Note that (\ref{aug2427}) implies that $0\le\alpha_{kN}\le 1$ and~$\alpha_{NN}=1$, so $\alpha_{kN}$
satisfy the assumptions that we needed in Section~\ref{sec:alpha}. 

On the other hand, given a random outcome  voting model of Section~\ref{sec:alpha}, with a collection of probabilities~$\alpha_{kN}$ that additionally have the property that the probabilities
$\alpha_{kN}$ are increasing in $k$, then we can obtain a random threshold model 
by setting
\be\label{jun2120}
\beta_{kN}=\alpha_{kN}-\alpha_{k-1,N}.
\ee
Note that
\be\label{jun2156}
\sum_{k=0}^N\beta_{kN}=\alpha_{NN}=1.
\ee

Monotonicity of $\alpha_{kN}$ in $k$ is a natural assumption, as it says that the larger number of children voted
$1$ the higher the probability that the parent votes $1$. Moreover, it is easy to see that given a polynomial nonlinearity $f(u)$
satisfying (\ref{aug2406}), we can always find
a collection of probabilities~$\alpha_{kN}$ that is increasing in $k$ and so that (\ref{mar1710}) holds. 
To see that, let us take a nonlinearity $f(u)$ that is a polynomial
of degree $N$ such that $f(0)=f(1)=0$. The construction of the
probabilities $\alpha_{kN}$ in the argument leading to Theorem~\ref{thm-aug2402} produces $\alpha_{kN}$ that are increasing
in~$k$ as long as the branching rate $\beta$ satisfies the condition
\be\label{aug2429}
\beta\ge {2}{N}\max_{k}|b_k[f]|,
\ee
that is slightly stronger than (\ref{aug2425}). This is because if $\alpha_{kN}$ are given by (\ref{aug2410}) and (\ref{aug2422}), then
\be\label{aug2428}
\bal
\alpha_{k+1,N}&=\farc{k+1}{N}+\mu_k=\farc{k+1}{N}+\frac{b_{k+1}[f]}{\beta}=\farc{1}{N} +\alpha_{k,N}+\frac{b_{k+1}[f]-b_k[f]}{\beta}
\\
&\ge \farc{1}{N} +\alpha_{k,N}-\farc{2}{\beta}\max_{k}|b_k[f]|\ge\alpha_{k,N}. 
\enbal
\ee
We have proved the following.
\begin{thm}\label{thm-aug2404}
Let $f(u)$ be a polynomial of degree $N$ such that $f(0)=f(1)=0$. Then, there exists a random threshold voting model representation  in terms of a purely
$N$-ary branching Brownian motion for the 
solution to the initial value problem (\ref{mar1708}) with the initial condition $g(x)$ that is continuous and satisfies~$0\le g(x)\le 1$
for all $x\in\Rm^d$. 
\end{thm}

\subsection{Recursive up the tree propagation models for other nonlinearities}\label{sec:recursive}

The random outcome and random threshold voting models apply to equations (\ref{aug402}) with nonlinearities $f(u)$ such that
\be\label{aug2602}
f(0)=f(1)=0.
\ee
The reason for this restriction is that the above assumption guarantees that if the initial condition~$u(0,x)=g(x)$ 
satisfies $0\le g(x)\le 1$ for all $x\in\Rm^d$, then
\be\label{aug2629}
\hbox{$0<u(t,x)<1$ for all $t>0$ and~$x\in\Rm^d$.}
\ee 
Thus, it is conceivable that 
$u(t,x)$ can be interpreted as a probability of some event. For a general polynomial $f(u)$ that does not satisfy (\ref{aug2602}), 
solutions to (\ref{aug2616}) do not necessarily satisfy (\ref{aug2629}), so there is no reason to expect that they can be interpreted as a probability.
However, we can 
replace the voting model interpretation by a recursive propagation up the genealogical tree $\cT(t)$ of the branching Brownian
motion that we now describe. 

Let 
\be\label{aug2604}
f(u)=f_0+f_1u+\dots+f_Nu^N
\ee
be a polynomial of degree $N$. Consider the corresponding symmetric polynomial of $N$ variables
\be\label{aug2606}
S_N(u_1,\dots,u_N)=f_0+\farc{f_1}{N}(u_1+\dots+u_N)+f_2\cnk{N}{2}^{-1}\sum_{k\neq j}u_ku_j+\dots+f_N \prod_{i=1}^Nu_i, 
\ee
so that 
\be\label{aug2610}
f(u)=S_N(u,\dots,u).
\ee
To build a solution to (\ref{aug402}) with $f(u)$ as above, we run a purely $N$-ary BBM, with an exponential clock running at the rate $\beta=1$, 
until a time~$t>0$. The particles $X_1,\dots,X_{N_t}$ that are present at the time $t$ are assigned the random values
$u_k=g(X_k(t))$. Here, $g(x)$ is a given continuous function. Then, we propagate the values up the genealogical tree $\cT(t)$ by assigning to each
parent the value
\be\label{aug2608}
u_{\rm parent}=S_N(u_1,\dots,u_N)+\farc{u_1+\dots+u_N}{N}.
\ee
Here, $u_1,\dots,u_N$ are the values that have been previously assigned to the $N$ children of the parent under consideration. This recursive procedure
allows us to define the value $u_{\rm orig}$ at the root of the tree, the original particle that was present at the time $t=0$ at the position $x\in\Rm^d$.
We set 
\be\label{aug2612}
u(t,x)=\E_x[u_{\rm orig}].
\ee
Once again, the renewal argument using the independence of the offspring particles, nearly identical to that 
in (\ref{mar1706})-(\ref{mar1708}), shows that $u(t,x)$ satisfies
the initial value problem
\be\label{aug2614}
\bal
&u_t=\Delta u+f(u),\\
&u(0,x)=g(x).
\enbal
\ee
This gives the following.
\begin{thm}\label{thm-aug2602}
Let $f(u)$ be a polynomial of degree $N$. Then, there exists a recursive up the tree propagation representation in terms of a purely
$N$-ary branching Brownian motion for the 
solution to the initial value problem (\ref{aug2614}) with the initial condition $g(x)$ that is continuous and  bounded. 
\end{thm}

We point out that this model is deterministic and so it is not a special case of either of the two random voting models above.

\section{Examples of voting models}\label{sec:examples}

In this section, we consider some examples of voting models, beyond the heat equation and the Allen-Cahn equation we have considered
above. 

\subsubsection*{Random outcome voting models for the McKean nonlinearities}

Let us go back to the McKean type nonlinearities as in (\ref{mar1631}). First, we note the elementary identity:
\be\label{aug2914}
	1 - (1-u)^n
		= (1-u + u)^n- (1-u)^n
		= \sum_{k=0}^{n}\cnk{n}{k}u^{n-k}(1-u)^k- (1-u)^n
		=\sum_{k=0}^{n-1}\cnk{n}{k}u^{n-k}(1-u)^k.
\ee
This allows us to write a nonlinearity of the form (\ref{mar1631}), using (\ref{aug2302}) and (\ref{aug2914}), as
\be\label{aug2402}
\bal
f(u)&=\beta\Big(1-u -\sum_{n=2}^N p_n(1-u)^n\Big)=\beta \sum_{n=2}^N p_n \Big(1-(1-u)^n-u\Big)
\\
&=\beta\sum_{n=2}^N p_n\Big(\sum_{k=0}^{n-1}\cnk{n}{k}u^{n-k}(1-u)^k-u\Big).
\enbal
\ee
Comparing to (\ref{mar1710}), we see that this corresponds to the random outcome voting model that is not really random: 
we have $\alpha_{0n}=0$ for all $n\ge 2$,
and 
\be\label{aug2433}
\alpha_{kn}=1,~~\hbox{ for all $1\le k\le n$.}
\ee
Therefore, the McKean nonlinearities come from a very simple voting rule: the parent
particle votes $1$ if and only if at least one of its children voted $1$.   This, of course, agrees with the familiar interpretation of the probability distribution
of the maximum of BBM in terms of the solution to the Fisher-KPP equation~\cite{JB,Bovier}.

\subsubsection*{Uniformly biased voting models}

Next, we introduce a uniform bias in the random outcome voting model (\ref{mar1716}) we have obtained for the standard heat equation.
A parent with 
$n$ children, out  of which~$k$ voted $1$, now votes $1$ with a ``uniformly biased'' probability  
\be\label{mar1720}
\alpha_{kn}=\farc{(1+\gamma)k}{n},~~0\le k\le n-1,~~\alpha_{nn}=1.
\ee
Here, $\gamma\ge 0$ is a parameter measuring the ``bias'' toward voting $1$ versus voting $0$. As we need to have $\alpha_{kn}\le 1$ for all $1\le k\le n-1$, 
the bias~$\gamma>0$
needs to satisfy
\be\label{mar1727bis}
\gamma\le \farc{n}{n-1}-1=\farc{1}{n-1}. 
\ee
In particular, if $\gamma$ is fixed, only finitely many $p_n$ may be non-zero. 
Using expression (\ref{mar1720}) for $\alpha_{kn}$
in~(\ref{mar1710}) and recalling (\ref{jan2814}) gives the corresponding nonlinearity as 
\be\label{mar1721}
\bal
f(u)&=\beta\sum_{n=2}^N p_n\Big(\sum_{k=0}^n\cnk{n}{k}\alpha_{kn}u^k(1-u)^{n-k} -u\Big)\\
&=\beta\sum_{n=2}^N p_n\Big(\sum_{k=0}^n\cnk{n}{k}\farc{(1+\gamma)k}{n}u^k -\gamma u^n-u\Big)=
\beta\gamma\sum_{n=2}^N p_n(u-u^n).
\enbal
\ee
Taking $\gamma=\beta^{-1}$ gives  
\be\label{mar1722}
f(u)=u-A(u),~~A(u)=\sum_{k=2}^N p_ku^k. 
\ee
As we have seen in Section~\ref{sec:mck-notkpp}, these nonlinearities are of the Fisher-KPP  type but do not have
a McKean representation. 
Instead, such nonlinearities come from voting models with a uniform bias toward voting $1$, 
as in (\ref{mar1720}). Note that they lead to convex functions $A(u)$ such that $\alpha(u)=A(u)/u$
is also convex, which is impossible for the McKean type.

\subsubsection*{Group voting models}

The nonlinearities of the form (\ref{mar1721})  
have the property that $f'(0)\neq 0$ except in the trivial cases~$\beta=0$ or~$\gamma=0$.
In order to obtain a voting model
representation for nonlinearities with derivatives that vanish at $u=0$,
it is convenient to consider  voting  with a ``group-based'' bias. 
Let us fix some $m>1$ and assume that branching can only happen into $n>m$ children; that is, $p_k=0$ for all $k\le m$. 
The voting scheme is as follows: if a parent has $n$ children, of which $k$  vote $1$, and $k<m$, then the parent
votes $1$ with the unbiased probability 
\be\label{jan2725}
\alpha_{n,m}^{(k)} =\farc{k}{n},~~\hbox{if $0\le k<m$},
\ee
as in (\ref{mar1716}). However, if $m\le k\le n-1$, so that one can choose a group of $m$ out of $n$ children that all voted~$1$, then 
the parent
votes $1$ with the biased probability
\be\label{mar1724}
\bal
&\alpha_{n,m}^{(k)}(\gamma)=\farc{k}{n}+\gamma \cnk{k}{m}\cnk{n}{m}^{-1},~~\hbox{if $m\le k<n-1$},
\enbal
\ee
and, finally, if all children voted $1$, then 
\be\label{mar1725}
\bal
&\alpha_{n,m}^{(k)}=1,~~\hbox{if $k=n$.}
\enbal
\ee
Thus, the bias  
for the parent to vote $1$ relative to the unbiased probability $k/n$ 
is proportional to the ratio of the number ${k \choose m}$ 
 of $m$-tuples such that all particles in the $m$-tuple voted $1$ to the 
total number~${n \choose m}$ 
of~$m$-tuples of the $n$ children. This is a generalization of the bias
in the voting model~(\ref{mar1720})
for the nonlinearity~$u-u^n$, where the~$m$-tuple 
is simply a single particle. This leads to the nonlinearity 
\be\label{mar1727}
\bal
f(u)&=\beta\sum_{n=m+1}^N p_n\Big(\sum_{k=0}^n\cnk{n}{k}\alpha_{nm}^{(k)}(\gamma)u^k(1-u)^{n-k} -u\Big)\\
&=\beta\gamma\sum_{n=m+1}^N p_n \sum_{k=m}^{n-1}\cnk{n}{k}\cnk{k}{m}\cnk{n}{m}^{-1}u^k(1-u)^{n-k}. 
\enbal
\ee
Here, we used (\ref{jan2814}) and (\ref{jan2725})-(\ref{mar1725}). We now proceed to simplify the right side of~(\ref{mar1727}).
Expanding the term $(1-u)^{n-k}$ gives
\be\label{mar1728}
\bal
	f(u)
		&=\beta\gamma\sum_{n=m+1}^N p_n \sum_{k=m}^{n-1}\cnk{n}{k}\cnk{k}{m}\cnk{n}{m}^{-1}u^k \Big(\sum_{q=k}^{n}\cnk{n-k}{q-k}(-1)^{q-k}u^{q-k}\Big) \\
		&=\beta\gamma\sum_{n=m+1}^N p_n \sum_{k=m}^{n-1}\sum_{q=k}^{n}\farc{n!}{k!(n-k)!}\farc{k!}{m!(k-m)!}\farc{m!(n-m)!}{n!}\farc{(n-k)!}{(q-k)!(n-q)!}
(-1)^{q-k}u^q\\
&=\beta\gamma\sum_{n=m+1}^N p_n \sum_{k=m}^{n-1}\sum_{q=k}^{n} \farc{(n-m)!}{(k-m)!(q-k)!(n-q)!}
(-1)^{q-k}u^{q}\\
&=\beta\gamma\sum_{n=m+1}^N p_n \sum_{k=m}^{n-1}\sum_{q=k}^{n-1} \farc{(n-m)!}{(k-m)!(q-k)!(n-q)!}
(-1)^{q-k}u^{q}\\
&\qquad+\beta\gamma\sum_{n=m+1}^N p_n \sum_{k=m}^{n-1}  \farc{(n-m)!}{(k-m)!(n-k)!}
(-1)^{n-k}u^{n} = f_1(u)+f_2(u). 
\enbal
\ee
In order to simplify this expression for $f(u)$, 
we use the identity\be\label{jan2624}
\sum_{k=0}^{n-1}\frac{(-1)^{k}n!}{(n-k)!~k!}
= (-1)^{n+1},
\ee
that holds for all $n\ge 1$ and can be obtained by expanding $(1-1)^n$. 
This allows us to write
\be\label{mar1729}
\bal
f_2(u)&=\beta\gamma\sum_{n=m+1}^N p_n \sum_{k=m}^{n-1}  \farc{(n-m)!}{(k-m)!(n-k)!}(-1)^{n-k}u^{n}\\
&=
\beta\gamma\sum_{n=m+1}^N p_n \sum_{k=0}^{n-1-m}  \farc{(n-m)!}{k!(n-k-m)!}(-1)^{n-k-m}u^{n}\\
&=
\beta\gamma\sum_{n=m+1}^N p_n(-1)^{n-m}(-1)^{n-m+1}u^n=-\beta\gamma\sum_{n=m+1}^N p_nu^n.
\enbal
\ee
For the first term in the right side of (\ref{mar1728}), we can write
\be\label{mar1730}
\bal
f_1(u)&=\beta\gamma\sum_{n=m+1}^N p_n \sum_{k=m}^{n-1}\sum_{q=k}^{n-1} \farc{(n-m)!}{(k-m)!(q-k)!(n-q)!}
(-1)^{q-k}u^{q}\\
&=\beta\gamma\sum_{n=m+1}^N p_n \sum_{q=m}^{n-1}\sum_{k=m}^{q} \farc{(n-m)!}{(k-m)!(q-k)!(n-q)!}
(-1)^{q-k}u^{q}\\
&=\beta\gamma\sum_{n=m+1}^N p_n \sum_{q=m}^{n-1}\farc{(n-m)!}{(n-q)!}(-1)^qu^q\sum_{k=m}^{q} \farc{1}{(k-m)!(q-k)! }(-1)^k\\
&=\beta\gamma\sum_{n=m+1}^N p_n \sum_{q=m}^{n-1}\farc{(n-m)!}{(n-q)!}(-1)^{q+m}u^q\sum_{\ell=0}^{q-m} \farc{1}{\ell!(q-\ell-m)! }(-1)^{\ell}.
 \enbal
\ee
The last sum is, up to a $(q-m)!$ factor, a binomial expansion:
\be\label{mar1731}
\bal
	\sum_{\ell=0}^{q-m} \farc{1}{\ell!(q-\ell-m)! }(-1)^{\ell}
		=\farc{1}{(q-m)!} (1 - 1)^{q-m}.
\enbal
\ee
Thus, the only nontrivial term in $f_1$ occurs when $q=m$, leading to
\be\label{mar1732}
\bal
f_1(u)&=\beta\gamma \sum_{n=m+1}^N p_n u^m=\beta\gamma u^m.
 \enbal
\ee
Combining (\ref{mar1728}), (\ref{mar1729}) and (\ref{mar1732}) gives
\be\label{mar1733}
f(u)=\beta\gamma\sum_{n=m+1}^N p_n(u^m-u^n).
\ee
Thus, the group voting models lead to this simple class of nonlinearities that vanish at least quadratically at the origin.

\subsubsection*{The pushmi-pullyu transition for the Ebert-van Saarlos nonlinearities}

Let us now describe a voting model that gives a probabilistic interpretation for the nonlinearities 
\be\label{feb118}
f(u)=(u-u^n)(1+\chi nu^{n-1})=u-u^{2n-1}+(\chi n-1)(u^n-u^{2n-1}),
\ee
introduced in~\cite{Ebert-vanSaarlos}. 
This paper used a formal matched asymptotic expansions argument to show that these nonlinearities are of the pulled type for $0\le \chi<1$
and of the pushed type for $\chi>1$.  This was recently proved in~\cite{AHR-RD}, where it was also shown that they exhibit the so-called 
pushmi-pullyu behavior at $\chi=1$. As we have mentioned,  equations of the form (\ref{aug402})
are pulled if the spreading of the solutions is dominated by the growth in the regions where the solution is small, and
in the pushed regime it is dominated by the regions where the solution is close neither to $0$ nor to $1$.
We refer to~\cite{AHR-RD} for a detailed discussion of the pushed and pulled regime and the transition between them.
In particular, Bramson's universal ``extremals of log-correlated fields'' asymptotics (\ref{aug2624}) for the front location holds for all
$0\le\chi<1$, while for $\chi>1$ the front location has the pushed type asymptotics
\be\label{aug2902}
X(t)=c_*t+x_0,~~\hbox{as $t\to+\infty$},
\ee
without a logarithmic in time correction. Here, the speed $c_*$ is given by
\be
c_*=\sqrt{\chi}+\farc{1}{\sqrt{\chi}}.
\ee
Let us now give an interpretation of the critical value $\chi=1$ where the transition from the pulled to the pushed regime 
happens, from the point of view of voting models.  The interest here is mainly to understand what kind of observables of the 
branching trees in other log correlated fields would have the universal extremal behavior  (\ref{aug2624}) and what kinds
should be considered as non-extremal or bulk quantities, with asymptotics as in (\ref{aug2902}). 

To get a voting model for (\ref{feb118}), consider a BBM with $2n-1$ children born at each branching event, with the exponential clock running at rate $1$.
In addition, each parent is assigned an index $I$ or $G$
with the probabilities
\be
p_I=\frac{1}{1+\beta},~~p_G=\farc{\beta}{1+\beta},
\ee
with some $\beta>0$ fixed. 
The voting is done as follows. If a parent has the label $I$ and $N$
is the number of the parent's children that voted $1$, then the parent's probability to vote~$1$ has a uniform bias $\gamma>0$, as in (\ref{mar1720}):
\be\label{feb116}
\Pm(V_{\rm parent}=1\vert N=k\big)=\farc{(1+\gamma)k}{2n-1},~~\hbox{ for $0\le k\le 2n-1$,}
\ee
and
\be
\Pm(V_{\rm parent}=1\vert N=2n-1\big)=1. 
\ee
This is the pulled, or ``extremal log-correlated,'' component of the voting.  

If the vertex has the label $G$, then the parent's vote is based on the group voting, as in~(\ref{mar1724}). 
More precisely, let $N$ be the number of children that voted
$1$, out of the $2n-1$ particles. 
If $0\le N\le  n-1$, then the parent votes~$1$ with the unbiased probability
\be
\Pm(V_{\rm parent}=1\vert N=k\big)=\farc{k}{2n-1},~~\hbox{ for $0\le k\le n-1$.}
\ee
If all offspring particles voted $1$ so that $N=2n-1$, then the parent votes $1$:
\be
\Pm(V_{\rm parent}=1\vert N=2n-1\big)=1. 
\ee
Finally, if~$n\le N\le 2n-2$, then 
the parent votes~$1$ according to the group voting rules for the nonlinearity~$u^n-u^{2n-1}$: 
\be
\Pm\big(V_{\rm parent}=1\vert N=k\big)=\farc{k}{2n-1}+\gamma\cnk{k}{n}\cnk{2n-1}{n}^{-1}.
\ee
This is the pushed, or ``bulk log-correlated," component of the voting. 
 
This leads to the nonlinearity 
\be
f(u)=\gamma\Big(\farc{1}{1+\beta}(u-u^{2n-1})+ \farc{\beta}{1+\beta} (u^n-u^{2n-1})\Big),
\ee
which is a multiple of (\ref{feb118}). As the pushed-pulled nature of nonlinearity does not change after multiplying  by a factor, the pushmi-pullyu
transition at $\chi=1$ corresponds to $\beta=n-1$, or 
\be\label{aug2908}
p_I=\farc{1}{n},~~p_G=\farc{n-1}{n}.
\ee
The above voting procedure can be generalized for other log-correlated fields. 
It would be interesting to understand if the same bulk to extremal behavior transition happens for observables of log-correlated fields
other than BBM.

\bibliographystyle{plain}

\end{document}